\documentclass[12pt]{article}

\usepackage{amsmath}
\usepackage{amsfonts}
\usepackage{amssymb}
\usepackage{amsthm}
\usepackage{latexsym}
\usepackage{color}
\usepackage[T1]{fontenc}
\title{On one variant of strongly nonlinear  Gagliardo-Nirenberg inequality involving Laplace operator with application
 to nonlinear elliptic problems}

\author{Tomasz Choczewski$^{\rm a}$
\ and Agnieszka Ka\l{}amajska$^{\rm a,b,}$\thanks{Corresponding author. Email: A.Kalamajska@mimuw.edu.pl
\vspace{6pt}}\hspace{5pt}$^{,}$\footnote{The work of A.K. was supported by the National Science Center (Poland), Grant 2014/
14/M/ST1/00600}\\
\vspace{-6pt}
\small{$^{a}${\em{Faculty of Mathematics, Informatics, and Mechanics, University of Warsaw, Warsaw, Poland}}};\\
\small{$^{\rm b}${\em{Institute of Mathematics, Polish Academy of Sciences, Warsaw, Poland}}}}

\renewcommand{\em}{\sl}

plus 6pt minus 12pt
\newcommand{\barint}{
         \rule[.036in]{.12in}{.009in}\kern-.16in
          \displaystyle\int  }
\def\r{{\mathbb{R}}}
\def\n{{\mathbb{N}}}
\def\rn{{\mathbb{R}^{n}}}

\begin{document}

\newtheorem{theo}{\bf Theorem}[section]
\newtheorem{coro}{\bf Corollary}[section]
\newtheorem{lem}{\bf Lemma}[section]
\newtheorem{rem}{\bf Remark}[section]
\newtheorem{defi}{\bf Definition}[section]
\newtheorem{ex}{\bf Example}[section]
\newtheorem{fact}{\bf Fact}[section]
\newtheorem{prop}{\bf Proposition}[section]
\newtheorem{prob}{\bf Problem}[section]


\makeatletter 
      \@addtoreset{equation}{section}
      \renewcommand{\theequation}{\thesection.\arabic{equation}}
\makeatother


\newcommand{\ds}{\displaystyle}
\newcommand{\ts}{\textstyle}
\newcommand{\ol}{\overline}
\newcommand{\wt}{\widetilde}
\newcommand{\ck}{{\cal K}}
\newcommand{\ve}{\varepsilon}
\newcommand{\vp}{\varphi}
\newcommand{\pa}{\partial}
\newcommand{\rp}{\mathbb{R}_+}
\newcommand{\hh}{\tilde{h}}
\newcommand{\HH}{\tilde{H}}
\newcommand{\ct}{{\cal T}}

\maketitle

\smallskip
  {\small Key words and phrases: Gagliardo-Nirenberg inequalities, interpolation inequalities,  Sobolev spaces, regularity, elliptic PDE's.\\

MSC (2000): Primary 46E35, Secondary 26D10
  }

  \begin{abstract}
\noindent
We obtain the inequality $$\int_{\Omega}|\nabla u(x)|^ph(u(x))dx\leq C(n,p)\int_{\Omega} \left( \sqrt{ |\Delta u(x)||{\cal T}_{h,C}(u(x))|}\right)^{p}h(u(x))dx,$$ where
$\Omega\subset \rn$ is a bounded Lipschitz domain,    $u\in W^{2,1}_{loc}(\Omega)$ is postive and obeys some additional assumptions,   $\Delta u$ is the Laplace operator,
${\cal T}_{h,C}(\cdot )$ is certain transformation of the continuous function $h(\cdot)$. We also explain how to apply  such inequality to deduce regularity  for solutions of nonlinear eigenvalue problems  of elliptic type for degenerated PDEs, with the illustration within the  model of electrostatic micromechanical systems (MEMS).
\end{abstract}

\section{Introduction}
In  \cite{ma}, Section 8.2.1 one finds the following inequality:
\begin{equation*}
\int_{\{ x: u(x)\neq 0\} }  \left(\frac{ |u^{'}(x)|}{u^{\frac{1}{2}}(x)}\right)^pdx \le \left( \frac{p-1}{|1-\frac{1}{2} p|}\right)^{\frac{p}{2}}
 \int_{\bf R} |u^{''}(x)|^{{p}}dx,
\end{equation*}
where $p> 2$, $u\ge 0$ is  smooth and compactly supported. 
Inspired by this inequality, the authors of \cite{akijp} have obtained more general inequality:
\begin{equation}\label{ex1}
\int_{\{x:u(x)\neq 0\}} \left(\frac{ |u^{'}(x)|}{|u(x)|^\theta}\right)^pdx \le \left( \frac{p-1}{|1-\theta p|}\right)^{\frac{p}{2}}
 \int_{\{x:u(x)\neq 0\}} \left(\frac{\sqrt{ |u(x)u^{''}(x)|}}{|u(x)|^\theta}\right)^{{p}}dx,
\end{equation}
where $2\le p <\infty,\theta \in \mathbf{R}$, under certain assumptions on $u$, which permit non-negative smooth compactly supported functions. 

In fact, inequality \eqref{ex1} is the special case of the more general inequality:
\begin{equation}\label{jedwym}
\int_{(a,b)}|u^{'}(x)|^ph(u(x))dx\leq C_p\int_{(a,b)}\sqrt{\left|u^{''}(x){\cal T}_h(u(x))\right|}^{p}h(u(x))dx,
\end{equation}
where  $-\infty\le a <b\le +\infty$, $u$   obeys some assumptions and we can permit $u\in C_0^\infty ((a,b)), u\ge 0$, $h(\cdot )$ is a given continuous  function,
${\cal T}_h(\cdot )$ is  certain  transform of $h(\cdot )$ (see Definition \ref{hh}), the constant $C_p$ does not depend on $u$.
When for example $h(t) = t^\alpha, \alpha >-1,$ then ${\cal T}_{h}(\lambda)=\frac{1}{\alpha +1}\lambda$, so it is proportional to $\lambda$ and in that case the inequality \eqref{jedwym} takes the form  \eqref{ex1}.

 After the substitution of $h\equiv 1$ in \eqref{jedwym} and the application of  H\"older's inequality,  we obtain  the classical Gagliardo-Nirenberg multiplicative inequality (\cite{ga,n1}):
\begin{eqnarray}\label{klasycznagn}
\left( \int_\r |u^{'}(x)|^{p} dx\right)^{\frac{2}{p}} \le C \left(  \int_\r |u(x)|^qdx\right)^{\frac{1}{q}} \left( \int_\r |u^{''}(x)|^rdx\right)^{\frac{1}{r}},\  {\rm whenever}\  \frac{2}{p}= \frac{1}{q}+ \frac{1}{r},
\end{eqnarray}
 where  $p\ge 2$. 
However, the classical inequality \eqref{klasycznagn} holds within all $1\le p\le\infty$.
We call inequality \eqref{jedwym} strongly nonlinear multiplicative inequality, because of its  nonlinear form and the link with the classical Gagliardo-Nirenberg inequality \eqref{klasycznagn}. 

\smallskip

Inequalities \eqref{ex1} and their more general variants \eqref{jedwym} were further developed in several directions.
In particular, regularity assumptions for the admissible function 
$u$ in \eqref{ex1} have been relaxed  in \cite{akikm}, while Orlicz variants 
of \eqref{jedwym}, having the form:
\begin{equation}\label{ogolne1} \int_\r M(|u^{'}(x)|h(u(x)))dx\leq  D\int_\r M\left( \sqrt{
|u^{''}(x){\cal T}_h(u)(x)| }\cdot h(u(x)) \right)dx, 
\end{equation}
where $M$ is convex, were obtained in \cite{akijp2}. In the special case of $M(t)=t^p$, one retrieves the equation \eqref{jedwym},  under the assumption $2\le p <\infty$.
It was possible to consider also to the case of $1<p<2$ in \cite{cfk}, at the cost of certain modification of inequality \eqref{ogolne1}. In particular, in place of 
the transformation ${\cal T}_h(u)$, one deals with some more general 
transformation, which is nonlocal in most cases.   

In this paper we are interested in the multidimensional variant of \eqref{jedwym}. 
Partial answer to that problem was presented in \cite{akitch}, where we obtained the inequality:
\begin{equation*}
\int_{\Omega}|\nabla u(x)|^ph(u(x))dx\leq C(n,p)\int_{\Omega} \left( \sqrt{ |\nabla^{(2)} u(x)||{\cal T}_{h,C}(u(x))|}\right)^{p}h(u(x))dx,
\end{equation*}
where  $C(n,p)=\left( p-1+\sqrt{n-1} \right)^{\frac{p}{2}}$, $\Omega\subseteq \rn$ has Lipschitz boundary and $n\ge 2$,  $u:\Omega\rightarrow \mathbf{R}^n$ is positive and  belongs to certain subset in the Sobolev space $W^{2,1}_{loc}(\Omega)$,  $\nabla^{(2)} u$ is the Hessian matrix of $u$ and
${\cal T}_{h,C}(u)$ is the transformation of the continuous function $h(\cdot)$ as in Definition \ref{hh}, under certain additional assumptions.
As kindly pointed to us by Patrizia Donato, it would be interesting to obtain the inequality having the form:
\begin{equation}\label{dimfreeommm}
 \int_{\Omega}|\nabla u(x)|^ph(u(x))dx  \leq A_\Omega \int_{\Omega } \left( \sqrt{ |P u (x){\cal T}_{h,C}(u(x))|}\right)^{p}h(u(x))dx,
\end{equation}
involving some elliptic operator $P$. In this paper we contribute to the  answer on that question, having in mind Laplace operator $Pu=\Delta u$,
see Theorems \ref{goal5} and \ref{goal6}.
Let us mention  that the `unweighted' ($h\equiv 1$) multidimensional  variant  of inequality \eqref{ogolne1}, the inequality: 
\begin{eqnarray*}
\int_\rn M(|\nabla u|)dx&\le& C\int_\rn M(|u||\nabla^{(2)}u|)dx,
 \end{eqnarray*} 
 was obtained earlier in \cite{akikpp}.

Most of the already discussed variants of inequality \eqref{ex1} have already been applied. In particular, the inequality \eqref{ogolne1} was applied in \cite{akijp2}
to  obtain second order Orlicz type izoperimetric inequalities and capacitary estimates.
Inequalities of the type \eqref{jedwym} were used to deduce regularity and asymptotic behaviour of solutions to general nonlinear eigenvalue problems for singular PDE's in \cite{akijp}.
Such an example  problem is the Emden-Fowler type equation
:
$
u^{''}(x)= g(x)u(x)^\alpha ,\  \alpha\in {\bf {R}}$, where we assume that $g(x)\in L^p((a,b))
,$ which appears for example in electricity theory, fluid dynamics or mathematical biology, see e. g. \cite{{NachmanCal},{carleman1957problemes},fermi,{luningperry2},{perry},tomas,{whitman},{wong75}}. Regularity of solutions to such a problem was studied in \cite{akikm}, by usage  of inequality  \eqref{ex1}. In another source, \cite{pesz}, the variant of  \eqref{ex1} was used  to
 obtain regularity of solutions of the   Cucker-Smale equation with  singular communication weight.
 
As explained in Section \ref{applications}, our new derived inequality: 
\begin{equation*}
 \int_{\Omega}|\nabla u(x)|^ph(u(x))dx  \leq A_\Omega \int_{\Omega } \left( \sqrt{ |\Delta u (x){\cal T}_{h,C}(u(x))|}\right)^{p}h(u(x))dx,
\end{equation*}
can be applied to deduce regularity for positive solutions to nonlinear PDE's like:
\begin{align*}
\left\{\begin{array}{ll}
\Delta u = g(x) \tau(u) & {\rm in} \quad \Omega, \\
u \equiv 0 & {\rm on} \quad \partial\Omega
\end{array} \right.
, 
\quad {\rm where}\ \ g\in L^{p/2}(\Omega), p\ge 2.
\end{align*}
In the particular case of $\tau (u) =\frac{1}{(1-u)^2}$, $0<u<1$, we deal with the 
 model of electrostatic micromechanical systems (MEMS). Using our proposed methods, we deduce that then the composition: $(1-u)^{\frac{1}{2}}$ belongs 
 to $W^{1,p}(\Omega)$ and moreover:
 \begin{equation*}
4\left( \int_\Omega |\nabla ((1-u)^{\frac{1}{2}})|^{p} dx\right)^{\frac{2}{p}} = 
\left( \int_\Omega |\nabla u|^{p} (1-u)^{-\frac{p}{2}} dx\right)^{\frac{2}{p}}\le (p-1) \left( \int_\Omega |g(x)|^{\frac{p}{2}}dx\right)^{\frac{2}{p}}.
\end{equation*}
We believe that the presented inequality  can be extended in many ways and applied 
to the regularity theory in nonlinear degenerated PDEs of various type.

\section{Preliminaries and notation}
{\bf Notation.}
From now on, we assume that $\Omega\subseteq \r^n$ is an open domain, $n\in \n$.
 In some cases $\Omega= B(0,r)$ is the ball of center at $0$ and radius $r$. When $\Omega$ has Lipschitz boundary, that is $\Omega$ belongs to the $\mathcal{C}^{0,1}$- class, see e.g. \cite{ma}, by $d\sigma$ we denote the $n-1$-dimensional Hausdorff measure on $\partial\Omega$. We use the standard notation: $C_0^\infty (\Omega)$ to denote smooth functions with compact support, $W^{m,p}(\Omega)$ and $W^{m,p}_{loc}(\Omega)$
to denote the  global and local Sobolev functions defined on $\Omega$, respectively. If $I\subseteq \mathbf{R}$ is an interval, by $AC(I)$ we denote functions which are absolutely continuous on $I$.
If $A\subseteq \r$ and $f$ is defined on $A$
by $f\chi_A$ we mean an extension of $f$ by zero outside set $A$. When $1<p<\infty$, we define  the continuous function  $\Phi_p:  \mathbf{R}^n\rightarrow \mathbf{R}^n$, by
\begin{align*}
\Phi_p(\lambda)=\left\{
\begin{array}{ccc}
|\lambda|^{p-2}\lambda & {\rm if}& \lambda\neq 0\\
0&{\rm if} & \lambda =0.
\end{array}
\right.
\end{align*}
When $a\in \r^n, b\in\r^n$, by $a\otimes b$ we denote the tensor product of $a$ and $b$, the square matrix  $(a_ib_j)_{i,j\in\{ 1,\dots ,n\} }$, while $a\cdot b$ denotes the scalar product of $a$ and $b$.
By $\nabla^{(2)} u(x)$   we mean Hessian matrix of
twice differentiable function $u:\Omega\rightarrow \rn$, the matrix $\left(\frac{\partial^{2}u (x)}{\partial x_i\partial x_j}\right)_{i,j\in \{ 1,\dots,n\} }$.
When $1<p<\infty$, then $p^{'}$ is H\"older conjugate to $p$, i.e. $p^{'}=p/(p-1)$. If $A$ is a vector or matrix, by $|A|$ we denote its Euclidean norm. 

We also recall the definition of the  infinity Laplacian: 
\begin{align*}
{\Delta_{\infty}}u(x):=\left\{
\begin{array}{ccc}
v(x)^t\nabla^{(2)}u(x)v(x), \ {\rm where}\  v(x)=\frac{\nabla u(x)}{|\nabla u(x)|} & {\rm if}& \nabla u(x)\neq 0\\
0 &{\rm if}& \nabla u(x)=0.
\end{array}
\right.
\end{align*}
We will consider certain operator involving $\Delta_\infty$-Laplacian:
\begin{align*}
\Delta^\spadesuit u (x) &:=&\Delta_\infty u(x)-\Delta u(x)~~~~~~~~~~~~~~~~~~~~~~~~~~~~~~~~~~~~~~~~~~~~~~~~~~~~~~~~~~~ \\&=&\sum_{i,j \in \{ 1,\dots,n\}, i\neq j}v_i(x)v_j(x)\frac{\partial^{(2)}u }{\partial x_i\partial x_j}(x)+
\sum_{i=1}^n ((v_i(x))^2-1)\frac{\partial^{(2)}u }{\partial x_i^2}(x),
\end{align*}
where $v(x)=\frac{\nabla u(x)}{|\nabla u(x)|}\chi_{\{ \nabla u(x) \neq 0\} }\in\mathbf{R}^n$.

\smallskip
\noindent
{\bf Transformation of weights.} In the sequel we will use the following definition, which has appeared first in  \cite{akijp}.

\noindent
\begin{defi}[transformations of the nonlinear weight $h$]\label{hh}
Let $0< B\le \infty$, $h:(0,B)\rightarrow (0,\infty)$ be a continuous  function which is  integrable on  $(0,\lambda)$ for every $\lambda <B$, $C\in\r$,
and let $H_C:[0,B)\rightarrow \mathbf{R}$ be the locally absolutely continuous primitive of $h$ extended to $0$, given by
\[
H_C(\lambda):= \int_0^\lambda h(s)ds-C,\ \ \lambda \in [0,B).
\]
We define
the transformation of $h$, ${\cal T}_{h,C}: (0,B)\rightarrow (0,\infty)$ by
\[
{\cal T}_{h,C}(\lambda):=
\frac{H_C(\lambda)}{h(\lambda)}, \ \ \lambda \in (0,B).
\]
\end{defi}
\noindent
In the case when  $C=0$ we omit it from the notation, i.e. we write $H_0=:H,{\cal T}_{h,0} =: {\cal T}_{h}$.
Note that $h$ and  ${\cal T}_{h,C}$ might not be defined at $0$ or $B$ in case $B<\infty$.

The following example shows that in many situations we may substitute the transformation ${\cal T}_{h}$  by $h$, at the cost of constant in the estimate.

\begin{ex}\rm
 When  $h(\lambda)=\lambda^\theta$, $\theta >-1$, then  ${\cal T}_{h}(\lambda)=(1+\theta)^{-1}\lambda$, so  ${\cal T}$ is proportional to the identity function. Similarly, ${\cal T}_{h}(\lambda)$  can be estimated from above
  by the proportional to the identity function,  when  we can deduce that  {${\cal T}_{h}(\lambda)= \frac{H(\lambda)}{h(\lambda)}\le A\lambda$,} with some general constant $A$. This is always the case when  $H$ is convex. 
\end{ex}

\smallskip
\noindent
{\bf The classical Gagliardo-Nirenberg interpolation inequality.} We recall the variant of the classical statement of the  Gagliardo and Nirenberg interpolation inequality (\cite{ga,n1}).

\begin{theo}\label{gnclass}
Let  $\Omega\subset \mathbf{R}^n$  be  the bounded Lipschitz domain,, $p,q,r\in [1,\infty
]$,
 $\frac{2}{q}
=\frac{1}{r} +\frac{1}{p}$,
 $0<k<m$ and $k,m$ are positive integers.
 Then there exists constant $C>0$ such that for any $u\in W^{2,1}_{\rm loc}(\Omega)$ 
\[
\|\nabla  u\|_{L^q(\Omega)}\le
C \left( \|u\|_{L^r(\Omega)}^{1/2}\|\nabla^{(2)}u\|_{L^p(\Omega)}^{1/2} + \|u\|_{L^r(\Omega)}\right). \]
\end{theo}

\smallskip
\noindent
{\bf First strongly nonlinear multidimensional interpolation inequalities.}
In \cite{akitch} we obtained the following variant of inequality  \eqref{jedwym}.  
We omit the formulation some more general statements proven there, also where we admitted functions $u$ such that $A<u<B$ with general constants $A,B$, or more general domains $\Omega$.



\begin{theo}\label{main2}
Let $n\ge 2$, $\Omega\subset \mathbf{R}^n$  be the bounded Lipschitz domain,
$2\le p<\infty$, $0<B\le \infty$, $h:(0,B)\rightarrow (0,\infty)$, $H_C, {\cal T}_{h,C}$  be as in Definition \ref{hh},
 $u\in  W^{2,p/2}(\Omega)\cap C(\bar{\Omega})$,
  $0<u<B$  in $\Omega$ and
  \begin{align}\label{aproksym1}
  \int_{\partial\Omega}\Phi_p(\nabla u(x))\cdot n(x) H_C(u(x))d\sigma (x)\in [-\infty,0],
  \end{align}
where $n(x)$ denotes unit outer normal vector to $\partial\Omega$, defined for $\sigma$ almost all $x\in \partial\Omega$.

Then
\begin{align*}
\int_{\Omega}|\nabla u(x)|^ph(u(x))dx\leq C(n,p)\int_{\Omega} \left( \sqrt{ |\nabla^{(2)} u(x)||{\cal T}_{h,C}(u(x))|}\right)^{p}h(u(x))dx,
\end{align*}
where  $C(n,p)=\left( p-1+\sqrt{n-1} \right)^{\frac{p}{2}}$.
\end{theo}

\begin{rem}\label{uwa}\rm  Assume that $H_C(0)\ge 0$.
When we deal with Dirichlet condition: $u=0$ on $\partial\Omega$, then $\nabla u$ is perpendicular to $\partial \Omega$. As  $u$ is positive inside $\Omega$, we have $\partial_{n(x)}u (x)=\nabla u(x)\cdot n(x) \le 0$ for $\sigma$ almost every $x\in\partial\Omega$. Therefore also $\Phi_p(\nabla u (x))\cdot n(x) \le 0$,
for $\sigma$ almost every $x\in\partial\Omega$, and the condition \eqref{aproksym1} is satisfied. 
\end{rem}

We also focus on the following statement from \cite{akitch}, the variant of Theorem 5,1, which was our inspiration for further work on the nonlinear interpolation inequalities.

\begin{theo}\label{goal3}
Let the assumptions of Theorem  \ref{main2} be satisfied.
Then we have
\begin{align}
\left( \int_{\Omega}|\nabla u(x)|^ph(u(x))dx \right)^{\frac{2}{p}} &\leq  (p-2)\left( \int_{\Omega} \left( \sqrt{ |({\Delta_{\infty}}u(x)){\cal T}_{h,C}(u(x))|}\right)^{p}h(u(x))dx\right)^{\frac{2}{p}} \nonumber\\& +  \left( \int_{\Omega} \left( \sqrt{ |\Delta u(x){\cal T}_{h,C}(u(x))|}\right)^{p}h(u(x))dx\right)^{\frac{2}{p}}  .\label{nier11}
\end{align}
In particular, when $p=2$, we have
\begin{align*}
 \int_{\Omega}|\nabla u(x)|^2h(u(x))dx  &\leq  \int_{\Omega}  |\Delta u(x)|{\cal T}_{h,C}(u(x))| h(u(x))dx .
\end{align*}
\end{theo}

\section{The question and one negative example}
It would be interesting to know if one could obtain stronger inequality than \eqref{nier11}, having the form:
\begin{align}\label{koercive}
 \int_{\Omega}|\nabla u(x)|^ph(u(x))dx \leq C_p \int_{\Omega} \left( \sqrt{ |\Delta u(x){\cal T}_{h,C}(u(x))|}\right)^{p}h(u(x))dx  ,
\end{align}
with some constant $C_p$ independent on  $u$, eventually positive inside $\Omega$ and such that $u\equiv 0$ on $\partial\Omega$, where $p>2$. The following example shows that such inequality is not possible when $\Omega$ is not bounded, even within the class of radially symmetric functions defined on the complement of the ball.

\begin{ex}\rm
 Assume that $p>2$, and let $\Omega = \rn \setminus B(0,1)$, $h(x) = x ^{-\alpha}$,  $\alpha <1$, in particular $H(x)= \frac{x^{-\alpha+1}}{1-\alpha}$. 
We consider
\begin{align*}
u(x) := \left\{
\begin{array}{ccc}
1 - \frac{1}{|x|^{n-2}} & {\rm when} & n \geq 3\\
\log |x| & {\rm when} & n=2
\end{array}
\right.
.
\end{align*}
 Obviously 
 $\Delta u =0$ in 
 $\Omega$, so that the right side in \eqref{koercive} is $0$, whereas the left one is positive. Moreover,
u satisfies the boundary condition \eqref{aproksym1}, because $u\equiv 0$ on $\partial\Omega$. 
\end{ex}

\noindent
Our next goal is  to give a positive  answer to our question within the class of radial functions defined on balls.

\section{Analysis within the radially symmetric functions}
In the whole section we assume that $\Omega$ is a proper ball with center at the origin. Moreover, we admit radial functions only to the inequality \eqref{koercive}.

 We start with the following statement, which gives the inequality slightly weaker than \eqref{koercive}. It contains the additional term depending on $u$ but not on its derivatives on the right hand side.

\begin{theo}\label{goal4}
Assume that $\Omega=B(0,r)\subseteq \r^n$,   $0<r,B<\infty$, $n\ge 2$, 
 $p> 2$, $C\in\mathbf{R}$,  let $h:(0,B)\rightarrow (0,\infty)$  and  $H_C$, ${\cal T}_{h,C}$ be  as in Definition \ref{hh}, $H_C(0)\ge 0$. Moreover, let 
\begin{align*}
u\in {\cal R}:=\{ u\in  W^{2,p/2}(\Omega)\cap C(\bar{\Omega}), :   u(x)=w(|x|),
0<w<B \ {\rm on}\ [0,r), w(r)=0
 \} 
\end{align*}
and 
\begin{align}
 \int_{\Omega } \frac{(\mathcal{T}_{h,C}(u))^p(h(u))^2}{|x|^p}dx<\infty.\label{er3}
 \end{align}
Then 
\begin{align}
\left( \int_{\Omega}|\nabla u(x)|^ph(u(x))dx \right)^{\frac{2}{p}} \leq  2 (p-1)\left( \int_{\Omega} \left( \sqrt{ |\Delta u (x){\cal T}_{h,C}(u(x))|}\right)^{p}h(u(x))dx\right)^{\frac{2}{p}} \nonumber\\+  [(p-2)(n-1)]^2 \left( \int_{\Omega} \left(  \frac{|{\cal T}_{h,C}(u(x))|}{|x|}\right)^{p}h(u(x))dx\right)^{\frac{2}{p}}. \label{coercive1}
\end{align}
\end{theo}
\noindent
{\bf Proof.} Let us denote: $u(x)=w(|x|)$.
Then we have
\begin{align*}
\nabla u(x) &= w^{'}(|x|)\cdot \frac{x}{|x|},\ |\nabla u(x)|=|w^{'}(x)|,\\
\nabla^{(2)} u(x) &= w^{''}(|x|) \frac{x}{|x|}\otimes \frac{x}{|x|} +\frac{w^{'}(|x|)}{|x|} \left[  I - \frac{x}{|x|}\otimes \frac{x}{|x|}\right],
\end{align*}
where by $I$ we denote the identity matrix. As $v^t (v\otimes v) v=1$ and $v^tIv=1$ when $|v|=1$, an easy computation gives:
\[
\Delta_{\infty}u(x)=w^{''}(|x|) \  {\rm and }\ \Delta u (x) = tr (\nabla^{(2)} u(x))= w^{''}(|x|) + (n-1)\frac{w^{'}(|x|)}{|x|} .
\]
As $u$ satisfies the assumptions of Theorem \ref{goal3} (see Remark \ref{uwa}), we have
\begin{align}
I&:= \left( \int_{\Omega}|\nabla u(x)|^ph(u(x))dx \right)^{\frac{2}{p}} 
\leq \nonumber\\& (p-2)\left( \int_{\Omega} \left( \sqrt{ |(w^{''}(|x|)){\cal T}_{h,C}(u(x))|}\right)^{p}h(u(x))dx\right)^{\frac{2}{p}}\nonumber\\& + \left( \int_{\Omega} \left( \sqrt{ |\Delta u(x){\cal T}_{h,C}(u(x))|}\right)^{p}h(u(x))dx\right)^{\frac{2}{p}}
=: I_1+I_2 .\label{pomoc}
\end{align}
Moreover,  $w^{''}(|x|) =\Delta u (x)- (n-1) \frac{w^{'}(|x|)}{|x|}$. Therefore 
\begin{align}
I_1&\leq (p-2)\left( \int_{\Omega} \left( \sqrt{ |\Delta u (x){\cal T}_{h,C}(u(x))|}\right)^{p}h(u(x))dx\right)^{\frac{2}{p}}\nonumber\\
&+ (p-2)(n-1)\left( \int_{\Omega} \left( \sqrt{ |\frac{|\nabla u(x)|}{|x|}{\cal T}_{h,C}(u(x))|}\right)^{p}h(u(x))dx\right)^{\frac{2}{p}}\label{pomoc2}\\ &\leq
(p-2)\left( \int_{\Omega} \left( \sqrt{ |\Delta u (x){\cal T}_{h,C}(u(x))|}\right)^{p}h(u(x))dx\right)^{\frac{2}{p}}\nonumber\\
 &+(p-2)(n-1)\left( \int_{\Omega}  {|\nabla u(x)|}^{p}h(u(x))dx\right)^{\frac{1}{p}}
\left( \int_{\Omega} \left(  \frac{|{\cal T}_{h,C}(u(x))|}{|x|}\right)^{p}h(u(x))dx\right)^{\frac{1}{p}}.\nonumber
\end{align}
Summing up $I_1$ and $I_2$ we obtain inequality:
\begin{align*}
I& \leq (p-1)\left( \int_{\Omega} \left( \sqrt{ |\Delta u (x){\cal T}_{h,C}(u(x))|}\right)^{p}h(u(x))dx\right)^{\frac{2}{p}} \\ &~~~~~~~~~~~~~~~~~~~~~~~~~~~~~~~~~~+ (p-2)(n-1)I^{\frac{1}{2}}
\left( \int_{\Omega} \left(  \frac{|{\cal T}_{h,C}(u(x))|}{|x|}\right)^{p}h(u(x))dx\right)^{\frac{1}{p}}.
\end{align*}
To finish the proof of \eqref{coercive1} we note that when $I<\infty$, then inequality $I\le a+ bI^{\frac{1}{2}}$ implies $I\le \left(\frac{b+ \sqrt{b^2+4a}}{2}\right)^2$, while $\left(\frac{b+ \sqrt{b^2+4a}}{2}\right)^2\le
 b^2 + 2a$.  
In the case of $I=\infty$ we have $I_1=\infty$ or $I_2=\infty$ in \eqref{pomoc}. When $I_2=\infty$, then \eqref{coercive1} holds trivially. 
The case of $I_2<\infty$ and $I_1=\infty$ is impossible 
because of \eqref{pomoc2}. Indeed, we then would have:
\begin{align*}
\infty=& \int_{\Omega } \left( \frac{|\nabla u(x)||\mathcal{T}_{h,C}(u(x))|}{|x|} \right)^{p/2} h(u(x))dx\\
 \le & \left( \int_{\Omega } |\nabla u(x)|^pdx     \right)^{1/2} \left( \int_{\Omega } \left( \frac{|\mathcal{T}_{h,C}(u(x))}{|x|} \right)^{p} (h(u(x)))^2dx \right)^{1/2} <\infty.
\end{align*}
The above expression is finite because of \eqref{er3} and an application  of the classical Gagliardo-Nirenberg's  inequality (Theorem \ref{gnclass}). Namely, when  $u\in L^\infty(\Omega)\cap W^{2.p/2}(\Omega)$, then $u\in W^{1,p}(\Omega)$. 
 \hfill$\Box$
 
\begin{rem}\rm ~\\
{\bf 1)} We do not know if technical assumption \eqref{er3}  can be violated in general, but in some cases it can be done.
Indeed, assume that 
\begin{align}\label{ghc}
G_{h,C}(\lambda):= \mathcal{T}_{h,C}(\lambda)h(\lambda)^{1/p}\in W^{1,\infty}((0,B) ).
\end{align}
 This implies that $v(x):= G_{h,C}(u(x))\in W^{1,1}_{loc}(\Omega)$  and $
\nabla v(x)= (G_{h,C}^{'}(u(x))\nabla u(x)\ {\rm a.e.}$ (see e.g. \cite{ma}, Section 1.1.3), 
so that
\[
 |\nabla v(x)|\le  \| G_{h,C}^{'}\|_{L^\infty(0,B)}|\nabla u(x)| \in L^p(\Omega).
\] 
Last conclusion follows from Theorem \ref{gnclass}, because $u\in L^\infty(\Omega),
\nabla^{(2)}u\in L^{p/2}(\Omega)$.  
The condition \eqref{ghc} is satisfied when we chose for example $h(\lambda)=\lambda^\alpha$ with $\alpha >0$ and $C=0$.\\
{\bf 2)} Another situation when \eqref{er3} holds is when $h\in L^\infty$. Indeed, let 
$[(p-2)(n-1)]^2\mathcal{A}^{\frac{2}{p}}$ be the second term in \eqref{coercive1}.
If $\mathcal{A}=\infty$ then \eqref{coercive1} holds trivially, while if $\mathcal{A}<\infty$, we have
$$
\int_\Omega \left( \frac{|\mathcal{T}_{h,C}(u(x))}{|x|} \right)^{p} (h(u(x)))^2dx \le 
 \| h\|_{\infty}\mathcal{A} <\infty .
$$
\end{rem}

\noindent
Our next statement allows to violate last term in \eqref{coercive1}, under some additional assumptions.

\begin{theo}
\label{goal5}
Let $\Omega=B(0,r)\subseteq \r^n$,  where $0<r,B< \infty$, 
 $2<p< n$, $C\in\mathbf{R}$, and $h:(0,B)\rightarrow (0,\infty)$,   $H_C$, ${\cal T}_{h,C}$ be  as in Definition \ref{hh}. 
Assume further that:\\
a) $H_C(0)\ge 0$,\\
b) $G_{h,C,p}(\lambda):={\cal T}_{h,C}(\lambda) h^{\frac{1}{p}}(\lambda): (0,B)\rightarrow (0,\infty)$ extends to locally absolutely continuous function defined on $[0,B)$ such that $G_{h,C,p}(0)=0$ and   
\begin{align*}
C_{h,C,p}:= {\rm sup}\{ \frac{|G_{h,C,p}^{'}(\lambda )|}{h^{\frac{1}{p}}(\lambda)} : \lambda\in (0,B)\} <\infty,\ 
D=D_{h,C,p,n}:=  (p-2)\frac{(n-1)p}{(n-p)}C_{h,C,p}  <1.
\end{align*}
  Moreover, let $u\in  W^{2,p/2}(\Omega)\cap C(\bar{\Omega})$ be the radially symmetric function such that   $0< u<B$ in $\Omega$, $u\equiv 0$ on $\partial\Omega$ and 
$
 \int_{\Omega } \frac{(\mathcal{T}_{h}(u))^p(h(u))^2}{|x|^p}dx<\infty\} .
$

\smallskip
Then  we have
\begin{align}\label{radialcoercive}
 \int_{\Omega}|\nabla u(x)|^ph(u(x))dx  \leq   A \int_{\Omega} \left( \sqrt{ |\Delta u (x){\cal T}_h(u(x))|}\right)^{p}h(u(x))dx,
\end{align}
where  $A=\left( \frac{2(p-1)}{1-D^2 }\right)^{\frac{p}{2}}$.
\end{theo}

\noindent
The proof will be based on the classical Hardy inequality (see e.g. \cite{hlp},
\cite{kmp}). It appears that best constants in the Hardy inequality have their impact on the choice of the admissable to the inequality 
\eqref{radialcoercive} set of functions $h$. More precisely, the term $\frac{(n-1)p}{n-p}$ in the definition of constant $D$ comes from the analysis of Hardy constant.

\begin{theo}\label{hardyclas}
Let $1<p<\infty$, $\alpha\neq p-1$. Suppose that $f=f(t)$ is locally
absolutely continuous function in $(0,\infty)$ such that
$\int_{0}^\infty |f'(t)|^p\,t^{\alpha}dt <\infty,$ and let
\begin{align}
f^+(0)&:= \lim_{t\to 0} f(t) =0\ \ {\rm for}\ \ \alpha <p-1,\nonumber\\
f(\infty)&:= \lim_{t\to\infty} f(t)=0\ \ {\rm for}\ \ \alpha
>p-1.\label{jekil}
\end{align}
Then the following inequality holds:
\begin{align*}
\int_{0}^\infty |f(t)|^pt^{\alpha -p}dt \le C\int_0^\infty
|f'(t)|^pt^{\alpha}dt,
\end{align*}
where $C=\left( \frac{p}{|\alpha -p +1| }\right)^p$ is best possible.
\end{theo}

\noindent
{\bf Proof of Theorem \ref{goal5}.} 
We can assume that right hand side in \eqref{radialcoercive} is finite.
We have $u(x):=w(|x|)$ and $u\in W^{2, 1}(B(0,r))$. It implies   that 
$w$ and $w^{'}$ are locally absolutely continuous on $(0,r)$. 
The function
\[
v(x):=G_{h,C,p}(u(x))= G_{h,C,p}(w(|x|))=: \tilde{v}(|x|)
\]
is radially symmetric. We have $\tilde{v}(r)=(G_{h,C,p}\circ w)(r)$, in particular 
$\tilde{v}^{'}(s)= G_{h,C,p}^{'}(w(s))\cdot w^{'}(s)$ and, as $w^{'}\in L^\infty_{loc}(0,r)$, 
$G_{h,C,p}^{'}\in L^1_{loc}(0,r)$, we get $\tilde{v}^{'}\in L^1_{loc}(0,r)$.
Therefore $\tilde{v}$ is also locally absolutely continuous on $(0,r)$. 
As $\tilde{v}(s)\to 0$ as $s\to r$, we can extend  $\tilde{v}$  by $0$ to the whole $(0,\infty)$, which gives locally absolutely continuous function. We keep the same notation for $\tilde{v}$ and for the extension.
Then $f=\tilde{v}$ satisfies condition (\ref{jekil}) with $\alpha =n-1$.  Applying polar coordinates and Theorem \ref{hardyclas}  with that $\alpha$
 we obtain:
\begin{align}\label{srodar}
{\cal A}&: =  \left( \int_{B(0,r)} \left(  \frac{|{\cal T}_{h,C}(u(x))|}{|x|}\right)^{p}h(u(x))dx\right)^{\frac{2}{p}} =
\left( \int_{B(0,r)} \left(  \frac{|v(x)|}{|x|}\right)^{p}dx\right)^{\frac{2}{p}}dx \nonumber\\ & = \left( \theta_{n} \int_0^\infty  | \frac{\tilde{v}(s)}{s}  |^p s^{n-1}ds\right)^{\frac{2}{p}} \le \theta_{n}^{\frac{2}{p}} \left( \frac{p}{n-p}  \right)^2\left(  \int_0^\infty |\tilde{v}^{'}(s)|^p s^{n-1}ds\right)^{\frac{2}{p}},
\end{align}
where $\theta_n$ is the $n-1$ dimensional measure of the unit sphere in $\mathbf{R}^n$.
Note that \[
|\tilde{v}^{'}(s)| =|G_{h,C,p}^{'}(w(s))w^{'}(s)|\le C_{h,C,p} h^{\frac{1}{p}}(w(s))|w^{'}(s)|,\] so that
\[
|\tilde{v}^{'}(s)|^p\le C_{h,C,p}^p|w^{'}(s)|^p h(w(s)).
\]
Using this and (\ref{srodar}), we obtain inequality:
\[
{\cal A}\le \left( \frac{p}{n-p}  \right)^2C_{h,C,p}^2\left( \int_{\Omega} |\nabla u(x)|^p h(u(x))dx\right)^{\frac{2}{p}}.
\]
According to Theorem \ref{goal4} this implies
\begin{align*}
I&:=\left( \int_{\Omega}|\nabla u(x)|^ph(u(x))dx \right)^{\frac{2}{p}} \\&\leq  2 (p-1)\left( \int_{\Omega} \left( \sqrt{ |\Delta u (x){\cal T}_h(u(x))|}\right)^{p}h(u(x))dx\right)^{\frac{2}{p}} + [(p-2)(n-1)]^2{\cal A} \\
&\leq  2 (p-1)\left( \int_{\Omega} \left( \sqrt{ |\Delta u (x){\cal T}_h(u(x))|}\right)^{p}h(u(x))dx\right)^{\frac{2}{p}} + D^2I=: \mathcal{X} +D^2I.
\end{align*}
By our assumption $\mathcal{A}<\infty$, therefore also $I<\infty$, because of first inequality above.
The condition $D<1$ implies the statement after we rearrange the inequality: $I\le 
\mathcal{X}+ D^2I$ with the finite $I$.  
 \hfill$\Box$

\noindent
The following remarks are in order.

\begin{rem}\rm ~\\
{\bf 1)} Let $h(\lambda)= \alpha \lambda^{\alpha -1}$, $\alpha >0$, $n>1$ and  $2<p< n$, $C=0$.  We have  
\begin{align*}
H(\lambda) =& \lambda^{\alpha}, \  h'(\lambda) = \alpha (\alpha -1) \lambda^{\alpha -2},\ \mathcal{T}_h(\lambda) = \alpha^{-1}\lambda,\ 
G_{h,0,p}(\lambda)= \alpha^{\frac{1}{p}-1}\lambda^{\frac{1}{p^{'}}+\frac{\alpha}{p}  },\\
G_{h,0,p}^{'}(\lambda)=&\alpha^{\frac{1}{p}-1}(\frac{1}{p^{'}}+\frac{\alpha}{p})
\lambda^{\frac{\alpha -1}{p}},\ \frac{G_{h,0,p}^{'}(\lambda)}{h^{\frac{1}{p}}(\lambda)}
= \frac{1}{\alpha}(\frac{1}{p^{'}}+\frac{\alpha}{ p})=: C_{h,0,p}.
\end{align*}
 Easy  computations show that conditions in Theorem \ref{goal5}  are satisfied when
\begin{align*}
D=(p-2)\frac{(n-1)p}{(n-p)}\frac{1}{\alpha}(\frac{1}{p^{'}} + \frac{\alpha}{ p})<1.
\end{align*}
We observe that $D <1$ if for instance when $p$ is sufficiently close to $2$.
\\
{\bf 2)}  When for some $E\in \mathbf{R}_+$ and  for any $\lambda\in (0,B)$ we have 
$
| H_C(\lambda)h^{'}(\lambda) |\le E  h^2(\lambda ),$ then
$C_{h,C,p}\le 1+(1-\frac{1}{p})E.$\\
{\bf 3)} In all cases $D<1$ for $p$ sufficiently close to $2$, provided that 
$C_{h,C,p}<\infty$. 
\end{rem}

\section{Relaxing the radiality assumption}
Our next results  apply to the case when $\Omega$ is bounded domain with $C^2$ boundary and $u$ is not necessarily radial.

\smallskip
\nonumber
We recall the following fact. For reader's convenience we enclose the sketch of the proof.

\begin{lem}\label{riesz1om}
Let $\Omega\subseteq \mathbf{R}^n$ be a bounded domain with $C^2$ boundary and $1< q<\infty$. 
 Then for any $u\in W^{2,q}(\Omega)$ such that $u\equiv 0$ on $\partial\Omega$ we have 
 \begin{eqnarray}
\left( \int_{\Omega}  |\nabla^{(2)} w (x)|^q dx\right)^{\frac{1}{q}} & \le  \tilde{C}_{q,\Omega} \left( \int_{\Omega} | \Delta w|^q dx\right)^{\frac{1}{q}},\label{secrieszom}\\
\left( \int_{\Omega} |\Delta^{\spadesuit}w|^q dx\right)^{\frac{1}{q}} & \le  \tilde{D}_{q,\Omega} \left( \int_{\Omega} |\Delta w|^q dx\right)^{\frac{1}{q}},\label{drugieom}
\end{eqnarray}
where $\tilde{D}_{q,\Omega}\le  \tilde{C}_{q,\Omega}+1$ and $\tilde{C}_{q,\Omega}>0$.
\end{lem}
\noindent
{\bf Proof.}
``\eqref{secrieszom}:'' Let $G(x,y)$ be Green's function for $\Omega$. Then
$u(x)=T(\Delta u) (x)$ where $Tv(x):=  \int_{\Omega} G(x,y)v(y)dy$.
It is known that  elliptic regularity theory 
that operator $\frac{\partial^2}{\partial x_i\partial x_j}T: L^q(\Omega)\rightarrow L^q(\Omega)$ is bounded, see e.g. \cite{cz}.\\
``\eqref{drugieom}:'' As $|\Delta_\infty u(x)|\le |\nabla^{(2)} u(x)|$, we obtain from \eqref{secrieszom}:
\begin{align*}
\left( \int_{\Omega} |\Delta^{\spadesuit}w|^q dx\right)^{\frac{1}{q}} &\le  \left( \int_{\Omega} |\Delta_\infty w|^q dx\right)^{\frac{1}{q}} + \left( \int_{\Omega} |\Delta w|^q  dx\right)^{\frac{1}{q}}\\ &\le 
\left( \int_{\Omega} |\nabla^{(2)}w|^q dx\right)^{\frac{1}{q}} + \left( \int_{\Omega} |\Delta w|^q  dx\right)^{\frac{1}{q}} \\ &\stackrel{\eqref{secrieszom}}
{\le} 
(\tilde{C}_{q,\Omega}+1) \left( \int_{\Omega} |\Delta w|^q  dx\right)^{\frac{1}{q}}.
\end{align*}
\hfill$\Box$

\noindent
Our next lemma deals with compositions.
\begin{lem}\label{zlozenia1}
Let $\Omega$ be a bounded domain with $C^2$ boundary,
 $0<B<\infty$, $1< q<\infty$ and  $W: (0,B )\rightarrow \mathbf{R}$ be such that 
$W\in C^1((0,B))\cap L^\infty ((0,B))$,  $w:= W^{'}$, 
   $\widetilde{W}:[0,B)\rightarrow \mathbf{R}$, $\widetilde{W}(\lambda):=\int_0^{\lambda} W(s)ds$ - be the primitive of $W$ such that $\widetilde{W}(0)=0$.

\noindent
Moreover, let   $v\in W^{2,q}(\Omega)\cap C(\bar{\Omega})$ satisfies:  $0<v(x)<B$ in $\Omega$.

Then the composition
$\widetilde{W}(v(x))$ belongs to $W^{2,q}_{loc}(\Omega)\cap C(\bar{\Omega})$,
$\Delta^{\spadesuit}\widetilde{W}(v)\in L^q(\Omega)$ and  
\begin{eqnarray}\label{rachunki}
\nabla (\widetilde{W}(v))  & = &  W(v)\nabla v,\\
\nabla^{(2)} (\widetilde{W}(v))  &= & w(v) \nabla v \otimes \nabla v + W(v)\nabla^{(2)} v,  \nonumber\\
\Delta_\infty (\widetilde{W}(v))  & = &  w(v)|\nabla v|^2 + W(v)\Delta_\infty v,\nonumber\\
\Delta (\widetilde{W}(v)) & = & w(v)|\nabla v|^2 + W(v)\Delta v,\nonumber\\
\Delta^{\spadesuit}(\widetilde{W}(v)) & = & W(v)\cdot \Delta^{\spadesuit} v,\nonumber
\end{eqnarray}
on the set of full measure in $\Omega$.
\end{lem}
{\bf Proof.} Properties  \eqref{rachunki} are obvious. 
Regularity property for $\widetilde{W}(v(x))$ follow from them easily once we are convinced that $\nabla^{(2)}\widetilde{W}(v)$ belongs to $L^q_{loc}(\Omega)$. To verify this
we use the fact that $v\in L^\infty$,   and so, by Theorem \ref{gnclass} we have $\nabla v\in L^{2q}(\Omega)$, $\nabla v\otimes \nabla v\in L^{q}(\Omega)$, while $w(v)$ and $W(v)$ are  locally bounded. \hfill$\Box$

Our main result in this section is the following.
\begin{theo}\label{goal6} 
Assume that $2<p<\infty$, $\Omega\subseteq \mathbf{R}^n$ is a bounded domain with $C^2$ boundary, $n\ge 2$,
 $0<B<\infty$ and the following conditions a),b),c) hold:

\begin{description}
 \item[a)] $h:(0,B)\rightarrow (0,\infty)$,   $H_C$, ${\cal T}_{h,C}$ are  as in Definition \ref{hh}, and additionally $h\in C^1((0,B))$,  $H_C(0)\ge 0$;
\item[b)] there exists some non-negative function $G\in L^\infty ((0,B))$ such that
$W:=G$  requires the assumptions in Lemma \ref{zlozenia1} and
for some positive constants $c_1, c_2$ we have 
\begin{equation}
c_1G \le |{\cal T}_{h,C}(\lambda)| h^{\frac{2}{p}}(\lambda) \le c_2G\ {\rm on}\ (0,B) 
\tag{5.4} 
\label{controlg}
\end{equation}
\item[c)] $\tilde{D}_{q,\Omega}$ is as in \eqref{drugieom}, 
\(
E =E_{h,p,C}:= {\rm sup} \{
\frac{|G^{'}(\lambda)|}{h(\lambda)^{\frac{2}{p}}} : \lambda \in (0,B) 
\}<\infty\)  and 
\begin{equation}\label{kappak}
 \kappa = \kappa_{h,C,p,\Omega}:= (p-2)c_2\tilde{D}_{p/2,\Omega}E<1.
 \tag{5.5}
\end{equation}
\end{description}
Then for every  function $u\in W^{2.p/2}_{loc}(\Omega)\cap C(\bar{\Omega})$ such that $0< u(x)<B$ in $\Omega$, $u\equiv 0$ on $\partial\Omega$  
\begin{equation}\label{dimfreeom}
 \int_{\Omega}|\nabla u(x)|^ph(u(x))dx  \leq A_\Omega \int_{\Omega } \left( \sqrt{ |\Delta u (x){\cal T}_{h,C}(u(x))|}\right)^{p}h(u(x))dx,
 \tag{5.6}
\end{equation}
where  $A_\Omega =\frac{(p-2)c_2c_1^{-1}\tilde{D}_{p/2,\Omega}E + (p-1)}{1-\kappa}$, provided that $I(u):= \int_{\Omega}|\nabla u(x)|^ph(u(x))dx$ is finite. Finitness of $I(u)$ is not needed when $G$ in \eqref{controlg} is  constant function.
\end{theo}

\noindent
{\bf Proof.}
Let us denote:
\begin{eqnarray*}
 J(u):= \int_{\rn} \left( \sqrt{ |\Delta u (x){\cal T}_{h,C}(u(x))|}\right)^{p}h(u(x))dx.
\end{eqnarray*}
From Theorem \ref{goal3}, Remark \ref{uwa} and the observation that 
$\Delta_\infty = \Delta^{\spadesuit} +\Delta$, we get:
\begin{eqnarray*}
(I(u))^{\frac{2}{p}}&\le & (p-2) \left(\int_{\Omega}\left( \sqrt{ |\Delta^\spadesuit u (x)\mathcal{T}_{h,C}( u(x))|}\right)^{p}h( u (x))dx \right)^{\frac{2}{p}}  \\
&+&
(p-2) \left(\int_{\Omega}\left( \sqrt{ |\Delta u (x)\mathcal{T}_{h,C}(u(x))|}\right)^{p}h(u (x))dx \right)^{\frac{2}{p}}
\\&+& \left(\int_{\Omega}\left( \sqrt{ |\Delta u (x)\mathcal{T}_{h,C}(u (x))|}\right)^{p}h( u(x))dx \right)^{\frac{2}{p}}\\
& =: &(p-2)(A(u ))^{\frac{2}{p}} +
(p-1)(J(u))^{\frac{2}{p}}.
\end{eqnarray*}
According  to Lemma \ref{zlozenia1}
with $W=G$,  we get: $(A(u))^{\frac{2}{p}}:=$
\begin{eqnarray*}
 \left(\int_{\Omega} |\Delta^\spadesuit u {{\cal T}_{h,C}(u) h^{\frac{2}{p}}}(u) |^{\frac{p}{2}}dx\right)^{\frac{2}{p}}\le 
\left(\int_{\Omega} |\Delta^\spadesuit u c_2{G}(u) |^{\frac{p}{2}}dx\right)^{\frac{2}{p}} = c_2
\left(\int_{\Omega} |\Delta^\spadesuit  \tilde{G}(u) |^{\frac{p}{2}}dx\right)^{\frac{2}{p}}.
\end{eqnarray*} 
We use Lemma \ref{riesz1om} to estimate:
\[
\left( \int_{\Omega} |\Delta^\spadesuit \tilde{G}(u)|^{\frac{p}{2}} dx\right)^{\frac{2}{p}}\le \tilde{D}_{\frac{p}{2},\Omega} \left( \int_{\Omega} |\Delta \tilde{G}(u)|^{\frac{p}{2}} dx\right)^{\frac{2}{p}}.
\]
According to our assumptions on $G$:
\begin{eqnarray*}
|\Delta \tilde{G}(u)|& =&
 |G^{'}(u)|\nabla u|^2 + G(u)\Delta u|\\
 &\le & 
E|\nabla u|^2(h(u))^{\frac{2}{p}} + c_1^{-1}|\Delta u||\mathcal{T}_{h,C}(u)|(h(u))^{\frac{2}{p}}.
\end{eqnarray*}
Therefore
\begin{eqnarray*}
(A(u))^{\frac{2}{p}}&\le & c_2\tilde{D}_{p/2,\Omega}E(I(u))^{\frac{2}{p}} + c_2c_1^{-1}\tilde{D}_{p/2,\Omega}(J(u))^{\frac{2}{p}},\\
(I(u))^{\frac{2}{p}}&\le & (p-2)c_2\tilde{D}_{p/2,\Omega} E(I(u))^{\frac{2}{p}} + (p-2)c_2c_1^{-1}\tilde{D}_{p/2,\Omega}(J(u))^{\frac{2}{p}} + (p-1)(J(u))^{\frac{2}{p}}.
\end{eqnarray*}
Therefore, if $\kappa =(p-2)c_2\tilde{D}_{p/2,\Omega}E<1$ and $J(u)<\infty$, we can rearrange the last inequality and obtain the statement, while if $E=0$, finiteness assumption 
for $J(u)$ is not needed. 
\hfill$\Box$

\noindent
Let us discuss our last statement.

\begin{rem}\rm~\\
{\bf 1)}
When $h$ is bounded on $(0,B)$ and $u\in W^{2.p/2}_{loc}(\Omega)\cap L^\infty (\Omega)$, then $I(u)<\infty$, what follows from Theorem \ref{gnclass}. \\
{\bf 2)} Simpler but barely possible situation is when $G(\lambda)= {\cal T}_{h,C}(\lambda) h^{\frac{2}{p}}(\lambda)$, so that  the assumption \eqref{controlg} 
 holds with $c_1=c_2=1$. In general, the derivatives of ${\cal T}_{h,C} h^{\frac{2}{p}}$ might be essentially larger than the derivatives of function $G$. \\
{\bf 3)} 
When ${\cal T}_{h,C}(\lambda) h^{\frac{2}{p}}(\lambda)\equiv C$ is a constant function, then the assumption $I(u)<\infty$ can be omitted. This is for example possible within the class of functions like $h(\lambda)=(B-\lambda)^\alpha$ for some $\alpha$'s, as  we will illustrate in our next section.\\
{\bf 4)}  We can omit the assumption $I(u)<\infty$ when $G\equiv 1$ in \eqref{controlg}, that is when $$0<c_1\le |({\cal T}_{h,C} h^{\frac{2}{p}})(\lambda)|\le c_2<\infty$$ in $(0,B)$.\\
{\bf 5)}
We do not know if in general it is possible to violate the assumption $I(u)<\infty$ from the statement.  The analysis dealing with one-dimensional similar inequality given in Remark 6.2 in  \cite{akijp} would suggest that it might not necessary be the case.\\
{\bf 6)} 
Note that the constant $A_\Omega$ in Theorem \ref{goal6} is dimension free.\\
{\bf 7)}
It appears important to know what is best constant $\tilde{D}_{q,\Omega}$
in the inequality \eqref{drugieom}, because it has its impact on the condition \eqref{kappak} for the  weights $h$ admissible to the inequality \eqref{dimfreeom}.\\
{\bf 8)} The condition \eqref{kappak} is satisfied when $E<\infty$ and $p$ is sufficiently close to $2$.
\end{rem}

\section{Applications to the elliptic eigenvalue problems}\label{applications}
{\bf The general approach to deduce regularity.}
In this section we use our results to prove the regularity of solutions. We illustrate our approach in an exemplary MEMS model from paper by  Esposito 
\cite{esposito}.

Let us consider an eigenvalue problem in the form:
\begin{align*}
\left\{\begin{array}{ll}
\Delta u = g(x) \tau(u) & {\rm in} \quad \Omega ,\\
u \in {\mathcal R}, & 
\end{array} \right.
\end{align*}
where $0<u<B$ in $\Omega$, 
 $\Omega$ is bounded domain with $C^2$ boundary, $\mathcal{R}$ is the set admissible to Theorem \ref{goal6}.

We use similar technics as in \cite{akijp}, Section 7.
Namely, assume that $g \in L^{q}(\Omega)$ where $q>1$. Formally we want to find some function $h$ satisfying the following inequality
\begin{equation}\label{regu}
|g(x)|^{q}= \left| \frac{\Delta u(x)}{\tau (u(x))} \right| ^{q} \ge c |{\cal T}_{h,C}(u(x)) \Delta u(x) | ^{\frac{p}{2}} h(u(x)),
\tag{6.1}
\end{equation}
where $q=p/2$, ${\cal T}_{h,C}$ is as in Definition \ref{hh}, $c>0$ is some constant. It can be done by  looking for positive function $h$ and constant $c>0$ such that for every $\lambda\in (0,B)$
\begin{eqnarray*}
|\tau (\lambda)|^{-q}\ge c|{\cal T}_{h}(\lambda)|^{q} h(\lambda) =c |H(\lambda)|^{q} h^{1-q}(\lambda).
\end{eqnarray*}
 If we find the function satisfying above equality and some additional conditions from Theorem \ref{goal6}, then by equalities \eqref{dimfreeom} and \eqref{regu}, we obtain
\begin{equation}\label{wupe}
|\nabla u |^{p} h(u) \in L^{1}(\Omega),\  \hbox{\rm equivalently}\ W_p(u)\in W^{1,p}(\Omega),\ {\rm where}\ W_p(\lambda) =\int_0^\lambda h^{\frac{1}{p}}(s)ds.
\tag{6.2}
\end{equation}

\smallskip
\noindent
{\bf Application to the model MEMS.}
In  \cite{esposito} there appears some simple model of electrostatic micromechanical systems (MEMS), which is reduced to the following problem
\begin{equation}\label{espomodel}
\left\{\begin{array}{ll}
\Delta u = \frac{r f(x)}{(1-u)^{2}} & {\rm in} \quad \Omega, \\
u= 0 & {\rm on} \quad \partial \Omega, \\
0<u<1 & {\rm in} \quad \Omega, 
\end{array} \right.
\tag{6.3}
\end{equation}
where $r >0$, $f \leq 0$, $u \in C^{1}( \overline{\Omega}) \cap W^{2,2}(\Omega)$, $\Omega$ is open and bounded, sufficiently regular. Let us possibly weaken the assumptions on $u$, requiring  that 
$u\in W^{2,q}(\Omega)\cap C(\bar{\Omega})$, $1<q<\infty$. We also do not need information about the sign of $f$.

Let $g(x)=r f(x)$ and assume that $f \in L^{q}(\Omega)$, $q=p/2$, $p>2$. We want to find function $h$ and constant $C>0$ such that
\begin{eqnarray*}\label{hszacowanie}
|\tau (\lambda)|^{-q} = (1-\lambda)^{2q}\ge c |H_C(\lambda)|^{q}h^{1-q}(\lambda).
\end{eqnarray*}
with some constant $c>0$.
Let us propose 
\[
h(\lambda):=(1-\lambda )^\alpha ,\ {\rm where}\ \alpha \in \mathbf{R}, C= \frac{1}{1+\alpha}.
\]
Clearly, such function is always integrable near zero. Moreover,  we have:
\begin{eqnarray*}
H_{\frac{1}{\alpha+1}}(\lambda)& =& -\frac{1}{\alpha +1} (1-\lambda)^{\alpha +1},\ H_{\frac{1}{\alpha+1}}(0)= -\frac{1}{\alpha +1} > 0, \ {\rm when}\ \alpha <-1, \\
\mathcal{T}_{h,\frac{1}{\alpha+1}}(\lambda)&=&-\frac{1}{\alpha +1}\frac{(1-\lambda)^{\alpha +1}}{(1-\lambda)^\alpha}= -\frac{1}{\alpha +1}(1-\lambda),\\
G_\alpha&:=&G_{h,\frac{1}{\alpha+1}, 2q}(\lambda)=-\frac{1}{\alpha +1}(1-\lambda)^{\frac{\alpha}{q}+1},\\
G_\alpha^{'}(\lambda) &=& \frac{1}{\alpha +1} (\frac{\alpha}{q}+1) (1-\lambda)^{\frac{\alpha}{q}} .
\end{eqnarray*}
We choose $\alpha =-q$, so that $G=G_\alpha$ is constant  and $E=0$ in Theorem \ref{goal6}. We leave it to the reader  to check that the remaining conditions, which are needed to apply Theorem \ref{goal6}, are satisfied. We arrive at the following regularity result:
\[
\left( \int_\Omega |\nabla u|^{2q} (1-u)^{-q} dx\right)^{\frac{1}{q}}\le (p-1)r \left( \int_\Omega |f(x)|^{q}dx\right)^{\frac{1}{q}}.
\]
In  particular, by \eqref{wupe},
$$(1-u)^{\frac{1}{2}}\in W^{1,2q}(\Omega),$$
which is not visible directly from the model
\eqref{espomodel}.

\smallskip
\noindent
{\bf Acknowledgement.} This research was partially done when T.Ch. was PHD student at Faculty of Mathematics, Informatics and Mechanics at the University of Warsaw. It was continued  when A.K. had research position at Institute of Mathematics of the Polish
Academy of Sciences at Warsaw in academic year 2018/2019. She would like to thank IM PAN for hospitality. We also  want to thank  Patrizja Donato  for drawing our attention to the interesting inequality  \eqref{dimfreeommm}.

\end{document}